\documentclass[a4paper,11.5pt,reqno]{amsart}

\usepackage{graphicx}
\usepackage[utf8]{inputenc}
\usepackage{listings}
\usepackage{color}

\usepackage{lineno}
\modulolinenumbers[5]

\usepackage{hyperref}

\numberwithin{equation}{section}

\renewcommand{\phi}{\varphi}

\begin{document}

\title[multivalues of the AGM reveal periodicities]
  {Computed multivalues of AGM reveal periodicities of inverse functions}

\author[F. Lamarche]{Fran\c{c}ois Lamarche}
\address{Teledyne LeCroy Inc, Chestnut Ridge, NY, USA}
\email{francois.lamarche@teledyne.com}

\author[H. Ruhland]{Helmut Ruhland}
\address{Santa F\'{e}, La Habana, Cuba}
\email{helmut.ruhland50@web.de}



\begin{abstract}
The article shows how two choices are possible whenever computing the geometric mean, and 
the repetition of this process can in general yield 2-to-the power N different values when 
the choices are compounded in the first N steps of evaluation of the arithmetic-geometric mean. 
This happens not only in the simple AGM involved in the computation of the complete elliptic  
integral of the first kind, but also in analogous methods for the computation of the complete 
and incomplete elliptic integrals of the first and second kind.
\end{abstract}

\date{\today}

\maketitle

    \section{INTRODUCTION}
    \label{section:intro}
    
	Computation of Elliptic integrals with Arithmetic-geometric means has been historically a very 
	powerful tool~\cite{gauss, jameson}, thanks to its quadratic convergence. Such a convergence was vital 
	before the advent of modern computers with practically unlimited precision available in 
	milliseconds. The duplication formula method represents a modern variation that allows a fast 
	evaluation~\cite{carlson}. However, in this article, we choose to focus on the computation 
	of the historic Elliptic Integrals of the First and Second kind, both Complete and Incomplete. 
	The Complete elliptic Integral of the Third kind is expressible as a function of incomplete 
	elliptic integrals of the the First and Second Kind~\cite{cylsphere}; and we will 
	leave out the computation of the Incomplete Elliptic integral of the Third kind as well as the 
	general case represented by the Carlson symmetric forms - the discussion of their multivaluedness 
	will be the subject of a future article.

	Fast processing is one advantage that modern compters provide; large amount of memory to store 
	multiple values is another. In this article, we focus on this aspect, using repeatedly the 
	multivaluedness of the square root function (the key to the geometric mean) to calculate 
	simultaneously multiple values at each iteration of the arithmetic-geometric mean.  
 
	Before all, we present the EAGM, which allows to calculate Incomplete and Complete Elliptic 
	Integrals of the first and second kind; then we describe the patterns of multivaluedness of the 
	quantities it produces:

	First, we look at the Complete Elliptic Integral of the first kind. Its multivalues span a set 
	that share the property of being quarter periods of the associated Elliptic functions - in a 
	very specific pattern.

	Secondly, we will focus on Incomplete Elliptic Integrals of the First kind, whose multivalues 
	possess a regularly repeating pattern with rectangular symmetry on the complex plane. This is 
	found to reflect the double-periodicity of its inverse function, the Jacobian elliptic function sn.

	Thirdly, we will look at multivalues of the complete Elliptic integrals of the Second Kind, which 
	make up a pattern related to that of the Complete Elliptic Integrals of the First kind via different 
	scaling of the real and imaginary part.

	And fourth, we will look at the multivalues of the Jacobi Zeta function, which represents the 
	non-trivial part of the Incomplete Elliptic function of the Second Kind. The pattern of 
	multivalues found here is very complicated, but a subset reveals an important symmetry of that 
	function, which is related to Weierstrass elliptic function (and to the essential elliptic 
	function~\cite{adlaj})

	After the conclusions, in appendices, we provide listings of the matlab functions
%
%
	and study the MAGM algorithm in the context of multivalues.

    \vspace{-0,2cm}
%
%
    \section{AN EXTENDED AGM, THE EAGM}
    \label{section:eagm}
	Here some formulas are given that define a quartet $(a_n, g_n, u_n, v_n)$ which we call
	EAGM for an obvious reason. The first 2 columns build the classical AGM. This quartet is the 
	base to calculate the incomplete elliptic integrals of the first and second kind.

	In each recursion step calculate the following two roots:
	\begin{equation}
	 r_n = \sqrt {a_n g_n}  \; \;  \; \; \;  \; \;
	 s_n = \sqrt { (u_n + v_n)^2 - (a_n - g_n)^2}/2 \label{TheEAGMRoots}
	\end{equation}
  
	With these 2 roots calculate in each step a row of the quartet with:
	\begin{equation}
 	 a_{n+1} = \frac{a_n + g_n}{2}  \; \;  \; \; \;  \; \;
 	 g_{n+1} = r_n   \; \;  \; \; \;  \; \;
 	 u_{n+1} = \frac{u_n + v_n}{2}   \; \;  \; \; \;  \; \;
 	 v_{n+1} = s_n       \label{FourEAGMRecursions}
	\end{equation}

	With $\Delta = \sqrt {1 - k^2 {\sin^2 (\phi})}$ set the initial row of the quartet:
	$$ (a_0, g_0, u_0, v_0) = (1, \sqrt {1 - k^2}, 1 / \sin (\phi), \Delta / \sin (\phi) ) $$

	With these initial values we get the following identity independent from $k$ and $\phi$ for all signs of the square roots and all recursion steps $n$:
	\begin{equation} a_n^2 - g_n^2 = u_n^2 - v_n^2 \label{Identitiy_aguv} \end{equation}

	Writing $x_\infty$ for  $\lim \limits_{n \rightarrow \infty} (x_n)$, we recognize additional properties:
	
	\begin{itemize}
	\item the sequences $a, g$ have a common limit, the well known $AGM (a_0, g_0) = a_\infty = g_\infty$
	\item the sequences $u, v$ have a common limit, $u_\infty = v_\infty$.
	\item for $\phi = \pi / 2$ the sequence $u, v$ is a copy of the sequence $a, g$ i.e. $u_n = a_n$ and $v_n = g_n$ 
	\end{itemize}

	These allow us to calculate the incomplete elliptic integrals of first kind:
	\begin{equation} F (\phi, k) = \sin^{-1} (a_\infty / u_\infty) / a_\infty \label{DefIncF} \end{equation}

	This quartet also defines a quadratically converging sum for the Jacobi Zeta function:
	\begin{equation} Z (\phi, k) = \sum_{n=0}^{\infty} 2^n (u_n - v_n) \sqrt {u_n^2 - a_n^2} / u_n \label{DefJacobiZeta} \end{equation}

	In the equations (\ref{TheEAGMRoots}), the sign of the chosen root is not specified. To be more specific, we write:

	\begin{equation}
	 r_n = \sigma_n \sqrt {a_n g_n}  \; \; \;\; \; \;
	 s_n = \delta_n \sqrt { (u_n + v_n)^2 - (a_n - g_n)^2} / 2 \label{SignedEAGMRoots}
	\end{equation}

	Here is the convention that lets us define the sign unambiguously: For $r_n$, if $\sigma_n$ is $+1$, then under this convention, the square root symbol refers to the one root such that 
	the division by $a_n+g_n$ gives a positive real part. So, for $\sigma_n=+1$, we obtain the root and geometric average nearer the arithmetic average, and for $\sigma_n=-1$, we obtain
          the root and geometric average father from the arithmetic average. When the two root are equally near/far, the positive imaginary part is chosen by convention. Note that taking the "far" 
          root in the calculation of the next AGM hinders convergence, while taking the "near" root in the calculation of the next AGM favors convergence. Note that, in the programs, a variable 
	called "sig g" takes care of flipping the sign in the programs if the real part of the sum is negative.

	For $s_n$, the convention is also geared on convergence. Under this convention, the square root symbol refers to the one of the two roots such that the addition of $u_n$ to this root
	makes the absolute value larger. Using this convention, $\delta_n$ being "+1"  will favor convergence toward a finite value, while $\delta_n$ being "-1" will hinder convergence to a 
	finite value and instead drive the common limit closer to zero. In the code, we achieve this convention by first taking the square root with a positive real part. Laying out the numbers 
	on the complex plane, if the next u value, $u_{n+1}$, is within 180 degrees, we don't need to change the sign. Otherwise, in the programs, by setting the variable called "sig" to "-1", 
           we change the sign.

	There is yet one more root that can be plus or minus: in (\ref{DefJacobiZeta}), we need to pick a convention for the square root sign. 
	By rewriting $u_n^2 - a_n^2$ as $(u_n - a_n)(u_n + a_n)$, we see that we are dealing with geometric mean of $(u_n - a_n)$ and $(u_n + a_n)$, so we can define the
           square root symbol as yielding the root near the arithmetic average $u_n$; then $\gamma_n$ defines if we put a negative sign, and so we represent all the possible choices by writing:
	\begin{equation} Z (\phi, k) = \sum_{n=0}^{\infty} 2^n (u_n - v_n) \gamma_n \sqrt {u_n^2 - a_n^2} / u_n \label{DefSignedJacobiZeta} \end{equation}

	\subsection{Analytic continuation and signs of square roots }

	In this article, we focus on the multivaluedness of complete/incomplete elliptic integrals of first and second kind. Having integrands $I_X (x)$ that are square roots of even 
	rational functions, they are 2-sheeted functions with 4 branch points.
	\begin{equation}
	X (z) =  \int \limits_{0}^{z} I_X (x) dx \quad \enspace I_X (x) = \sqrt {\frac{P (x)}{Q (x)}} \quad \enspace P (x), Q (x) \enspace \hbox{are even quadratic/quartics} \nonumber 
	\end{equation}
	Integration along different closed paths in the plane of the integration variable $x$ can yield different
	function values at a certain point $z$ in this plane and so result in the multivaluedness.
	Because the integrand $I_X (x)$ is 2-valued we have 2 different types of closed paths in the integration plane: 
	\begin{itemize}
	\item Type 1: the integrand takes the same values at the start and end points of the path 
	\item Type 2: the integrand takes different values (differing by a factor of $-1$) at the start and end points of the path, in other words, contrary to type 1, the path is not closed 
	in the 2-fold cover of the integration plane
	\end{itemize}

	The two different types of closed paths lead to two different behaviors of the analytically continued function at a point $z$: 
	\begin{itemize}
	\item for type 1, $X (z) \rightarrow + X (z) + c_1$ 
	\item for type 2: $X (z) \rightarrow - X (z) + c_2$ 
	\end{itemize}
	The constants $c_1, c_2$ do not depend on $z$.

	When a multivalued function is given by the limit of a recursion as in the case of $F (k)$ by the $AGM$
	or in the case of the Jacobi Zeta function (\ref{DefJacobiZeta}) by the recursion (\ref{TheEAGMRoots}) and (\ref{FourEAGMRecursions}), we can expect that analytic continuation 
	along 	different paths leads to different signs of the roots in the recursion. But not every sign combination has to belong to an analytic continuation along a closed path of the function.

	Because the integrand is branched at 4 points, there are two branch cuts, and a type 1 closed path can yield two different values $c_1$. This causes the multivalues of $X (z)$ 
	to belong to a 2-dimensional lattice, with the generating vectors being the two different values of $c_1$. In the case of the Jacobi Zeta function, one of the closed paths yields zero,
	and we get only a 1-dimensional lattice of multivalues.



	It is conjectured that for all signs of the roots in $r_n$, $s_n$ and for all multivalues of $\sin^{-1} ()$ in (\ref{DefIncF}) the corresponding function value can be reached 
	by analytic continuation, and that which  multivalue of  $\sin^{-1} ()$ we get corresponds to the type of contour:
	\begin{itemize}
	\item for type 1: $( + \sin_p^{-1} (a_\infty / u_\infty) +   n \pi) / a_\infty, n \enspace \hbox{even}$ 
	\item for type 2: $( -  \sin_p^{-1} (a_\infty / u_\infty) +   m \pi) / a_\infty, m \enspace \hbox{odd}$ 
	\end{itemize}
	where $\sin_p^{-1} ()$ is the principal value of the arc sine function.


    \vspace{-0,2cm}
%
%
    \section{THE INVERSE OF AGM MAKES A LATTICE ON THE COMPLEX PLANE}
    \label{section:K}
	We launch the computation of the AGM of unity with the parameter b (also known as k*) 
	using the MATLAB program "eagm.m" called by the program "fillk2x32.m" (see appendix
	\ref{appendix:program}). This MATLAB program, as well as its Excel 
	and python versions, which are attached as documents in the online version, implement a sign 
	choice $\sigma_i$, which can be $+1$ or $-1$, at each of the 5 first iterations, because the numbers 0 thru 31 are represented in binary by just 5 non-zero bits. The sign choices 
	$\sigma_i$ for $i>5$ are simply always "+1" in order to let convergence occur.

\begin{figure}
	\includegraphics[width=0.92\textwidth]{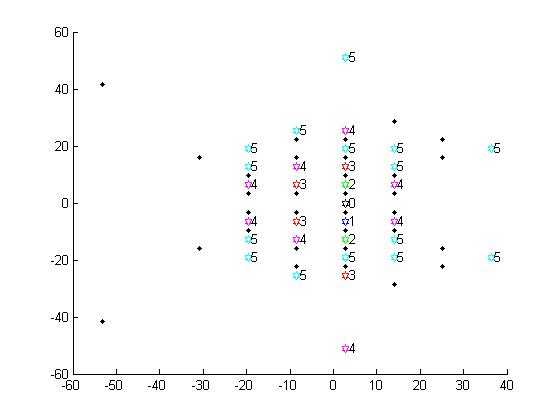} 
   	\caption{First thirty-two multivalues of $\pi/2/AGM$ on complex plane, for $b=\frac{1}{4}$, and for $b=-\frac{1}{4}$} \label{invagm}
    	\vspace{-0.2cm}
	\end{figure}

	As this program runs, it produces values for  $\pi/2/AGM$, shown by colored stars in the complex plane map~\ref{invagm}, which we notice are arranged in a 
	rectangular lattice surrounding and including K(k). Thinking in terms of to a "game of life"~\cite{conway} analogy, the points are labelled by the number of iterations of
	non-trivial sign choices needed to reach them, with many new points  reached each time the number of free sign choices is incremented by 1 
	(we stop at '5' here, and one '5' point at K(k)-64iK(b) is not shown as it is out of the plotting area).

	Also included in the figure are black points. They too form part of a rectangular lattice, but one that is interleaved vertically between the lattice of colored stars. They are 
	obtained when we start the AGM with unity and $g_0=-\sqrt{1-k^2}$ instead of $g_0=\sqrt{1-k^2}$.

    \subsection{Towards proofs for simple cases: an important relationship}

	The Landen transformation implies a relationship between $K(k)$ and $K(q)$ where $k=2\sqrt{q}/(1+q)$:
	$$K(k)=2/(1+b) K(q)$$
	with $b=\sqrt{1-k^2}$.
	This can be combined with the Landen transformation on the complimentary argument $v=\sqrt{1-q^2}$:
	$$K(v)=2/(1+q) K(b)$$
	to yield this relationship between the values of elliptic integrals:
	$$K(k)K(v)=2K(b)K(q)$$
	This result is because $(1+b)(1+q)$ is equal to two, as can be seen by expanding 
	$q=\sqrt{1-v^2}=\sqrt{1-(2\sqrt{b}/(1+b))^2}=(1+b)^{-1}\sqrt{(1+b)^2-4b}=(1+b)^{-1}(1-b)$, 
	which means that $1+q=(1+b)^{-1}(1+b+1-b)=2/(1+b)$. 
	Note that such a relationship can be extended over a second transformation; writing $q=2\sqrt{c}/(1+c)$ and $w=2\sqrt{v}/(1+v)$, we have:
	$$K(k)=\frac{2}{1+b}\frac{2}{1+v}K(c)$$
	and
	$$K(w)=\frac{2}{1+c}\frac{2}{1+q}K(b)$$
	and so 
	$$K(k)K(w)=4K(b)K(c)$$

	With this, we are now ready to perform an exact calculation of the effect of different sign choices on K(k). But first, we must note that there exist two ways to assess the 
	multivaluedness of $K(k)$ as calculated via $AGM(1,b)$. The first way is to consider that both b and -b are solutions to $k^2+b^2=1$. 
	This multiplies our overall sign choices by two, and, not surprisingly, it creates a denser lattice, with lattice parameters 4K(k) and 2iK(b).

    \subsection{A first way to calculate multivalued K(k)}
	For this, we just look at what $K(k)=\pi/2 / AGM(1,b)$ becomes if we choose $b_-=-\sqrt{1-k^2}$ for b. The denominator becomes
	$$AGM\left(1,b_-\right) = AGM\left( \frac{1+b_-}{2},\sqrt{b_-}\right) = \frac{1-b}{2} AGM\left(1, 2i\frac{\sqrt{b}}{1-b}\right)$$
	we can then use the definition of K to write it as:
	$$AGM(1,b_-) = \frac{1-b}{2} \frac{\pi}{2} \Big/ K\left( \sqrt{1-\left(2i\frac{\sqrt{b}}{1-b}\right)^2}\right) = \frac{1-b}{2} \frac{\pi}{2} \Big/ K\left(\sqrt{1+4b/(1-b)^2)}\right) $$
	$$ = \frac{1-b}{2} \frac{\pi}{2} \Big/ K\left(\frac{1+b}{1-b}\right) $$
	Known formulae~\cite{dlmf} allow us to compute K for k values greater than 1; these yield:
	$$AGM(1,b_-) = \frac{1+b}{2}\frac{\pi}{2} \Big/ \left(K\left(\frac{1-b}{1+b}\right) \mp i K\left(\sqrt{1-\left(\frac{1-b}{1+b}\right)^2} \right)\right)$$
	Using the variables introduced for the transformed value:
	$$AGM(1,b_-) = \frac{1+b}{2}\frac{\pi}{2} \Big/ (K(q) \mp i K(v))$$
	Finally, the ratio $\pi/2/AGM(1,b_-)$ can be finally written as:
	$$ \frac{2}{1+b}K(q) \mp 2i  \frac{1+q}{2} K(v)  = K(k) \mp 2i K(b)$$
	The distance between $K(k)$ and this value is $2iK(b)$, showing that the smallest complex numbers of the second lattice (-b instead of b) are offset by $2iK(b)$. 

    \subsection{A second, stricter way to calculate multivalued K(k)}
 
	We look at what happens if we decide the sign of b is not changeable.
	The denominator can still be called AGM(1,b), but in the calculation of the AGM, we choose the negative sign instead of the positive sign when calculating some of the geometric means. 			We denote this by $AGM_{-,+,...}$, where we indicate the chosen sign at each iteration, and the ellipsis (...) denotes a indefinite repetition of the last sign choice.
	As the computer algorithm runs, we explore many such sign choices, but in here, we will just explore what happens with only the first sign flipped; in this case, the denominator is
	$$AGM_{-,+,...}(1,b) = AGM\left(\frac{1+b}{2},-\sqrt{b}\right)= \frac{1+b-2\sqrt{b}}{4} AGM\left(1,\frac{\sqrt{(b+1)\sqrt{b}}}{1+b-2\sqrt{b}}i\right) $$

	This is again a call for an AGM with an imaginary number, corresponding to an evaluation of the complete elliptic integral of the first kind with an argument greater than one, 
	and so we again use the formula~\cite{dlmf}: $$AGM_{-,+,...}(1,b) = \frac{1+b}{2} \frac{1-v}{2} \frac{\pi}{2} \Big/ K\left(\frac{1+v}{1-v}\right) $$
	is
	$$\frac{1+b}{2}\frac{1+v}{2}\frac{\pi}{2} \Big/ \left(K\left(\frac{1-v}{1+v} \right) \mp i K\left(\sqrt{1-\left(\frac{1-v}{1+v}\right)^2} \right)  \right)$$
	and that is
	$$AGM_{-,+,...}(1,b) = \frac{1+b}{2} \frac{1+v}{2} \frac{\pi}{2} \Big/ (K(c) \mp i K(w)) = \frac{\pi}{2} \Big/ (K(k) \mp 4 i K(b)) $$
	What happens is the factor needed for returning K(w) to K(b) is just the same as the factor needed to return K(w) into K(b), except for the factor of four as shown above. 
	The difference between the scaled inverses of $AGM(1,b)$ and $AGM_{-,+,...}(1,b)$ is $4 i K(b)$, showing an example of the second lattice parameter being $4 i K(b)$ 
	(which is also the second period of all elliptic functions of parameter $k$).

    \subsection{Important observations and conjectures about the values obtained by a finite number of sign choices}

	In general, the choice of a square root of a product of two values can be divided between a square root which is nearer to the two values and one which is farther. 
	The exception is when the two values are exactly 180 degrees apart, i.e. opposite sign. In that case, one of the two roots (specifically, the root with a positive imaginary part)
	 is arbitrarily chosen as the nearest; for this reason, the pattern of stars does not have perfect up-down symmetry.

	The population of values obtained by considering more and more applications of different sign choices appear to constitute sort of a "game of life"~\cite{conway}, albeit 
	with very different rules.
	In contrast with the classical game, at each generation, the previously filled lattice points remain present. At each step, there are new multivalues, and they could in principle
	double the number of multivalues at each step. The population of multivalues will exactly double at each iteration of our zero-player game if there are no repeat values. 
	So far, it seems no point of the lattice exists repeatedly among the $2^n$ sign combinations, but we do not yet have a proof of that conjecture.

	The pattern of multivalues with up to N (in our example 5) arbitrary sign choices seems to produce values up to $2^{N-1}$ (in our example, $2^4=16$) unit lattice steps away. 
	The extent of the cloud of multivalues thus grows exponentially; however, it seems that some values on this lattice are never reached, such as non-trivial lattice points with 
	zero imaginary part. It is not yet known if all the non-zero imaginary part lattice points can be reached given a sufficient number of iterations of the sign choices.

	 It is speculated that a formal proof of the multivalue being exactly a lattice point could be developed automatically by a symbolic calculation generator, but we have not yet 
	attempted to develop such an automated proof builder.

	The values of $\frac{\pi}{2 AGM}$ reached by positive and negative values of b are very different, even though they both correspond to the same 
	$k=\sqrt{1-b^2}$ argument of the elliptic integral $K$. This suggests suggests that, in a sense, the function $\kappa (b)= \frac{\pi}{2 AGM(1,b)}$ is more fundamental 
	than $K(k)$. 

	The observations and conjectures about the numbered symbols in our graphic also apply to the black dots (where we start with -b); except those multivalues appears to have 
	exact up-down symmetry.

   \vspace{-2mm}
%
%
   \section{MULTIVALUE LATTICE FOR THE INCOMPLETE INTEGRAL OF THE FIRST KIND}
   \label{section:F}
 	We again launch the computation of the EAGM of unity with the parameter b (also known as k*) 
	using the MATLAB program "eagm.m" called by the program "fillf8x16.m" (see appendix
	\ref{appendix:program}). The EAGM has a second parameter, $\phi$, that controls the amount of incompleteness.
	To illustrate the working of this multivalue-computing program, we choose a value of b of 0.25 and a value of $\sin{\phi}$ of 0.8. We allow the sign of the geometric average $\sigma_n$ 
           to be +(the "near" geometric average) or - (the "far" geometric average) during the first 3  iterations; it is always chosen "+" after that, yielding convergence to a non-zero value. 
           We further allow the sign $\delta_n$ to be + or - during the first 4  iterations, and for the remaining  iterations it is always "+".

\begin{figure}
	\includegraphics[width=0.92\textwidth]{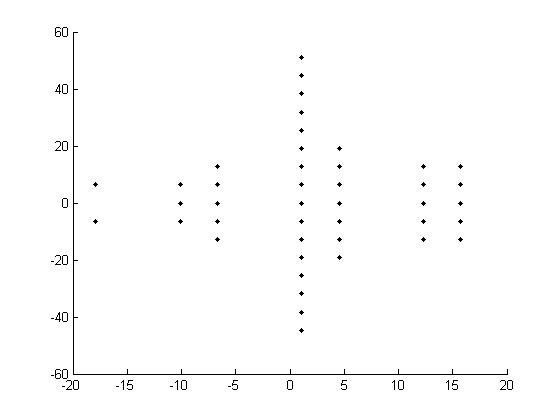} 
   	\caption{First 128 (8  by 16) multivalues of multivalued $F(\phi,k)$ on complex plane, for b=$\frac{1}{4}$ and $\sin{\phi}=0.8$} \label{eagm_f}
    	\vspace{-0.2cm}
	\end{figure}

	As this program runs, it produces values for  $\arcsin{z}/AGM$, shown by blue stars in the complex plane map~\ref{eagm_f}, which we can notice are arranged in a rectangular
           lattice with the same generating vectors as K(k).

          We observe that this lattice pattern comprises two multi-values per cell.

          We observe that sometimes, the same multivalues can be reached by different choices of $\sigma_n$ and $\delta_n$.

          We make a conjecture that every value on the lattice will be reached by sign combinations having a finite number of "-" choices.

	Even if not all values on this lattice of multivalues could be reached within a finite number of "-" sign choices, the $\arcsin{z}$ function has branches expressible as $\pi-\arcsin{z}$,
	as $2\pi+\arcsin{z}$, etc... and those will definitely fill the lattice. The reason is the multivaluedness of $K=\pi/2/AGM$ having the same lattice generating vectors as those of $F$.

          The lattice of values reflects the double periodicity of the inverse function of F: namely, the double-periodicity of the Jacobian elliptic function sn(z). But while the lattice generating values of 
	our F multivalues are $4K(k)$ and $4iK(b)$, the periods of sn are $4K(k)$ and $2iK(b)$. However, when considering all functions together, the Jacobi elliptic functions have periods 
	$4K(k)$ and $4iK(b)$.

   \vspace{-2mm}
%
%
   \section{MULTIVALUE LATTICE FOR THE COMPLETE INTEGRAL OF THE SECOND KIND}
   \label{section:E}
   
	Following~\cite{jameson}, the computation of the AGM can be accompanied by a summation of the square differences between the arithmetic and geometric means. The 
	square differences factor in the sum is not constant, but raises by a factor of two at each iteration.

	These sums are then divided by the AGM. This process is repeated in much the same way as the computation of the inverses of the AGM over choices of the sign of the geometric mean.

	\begin{figure}
	\includegraphics[width=0.92\textwidth]{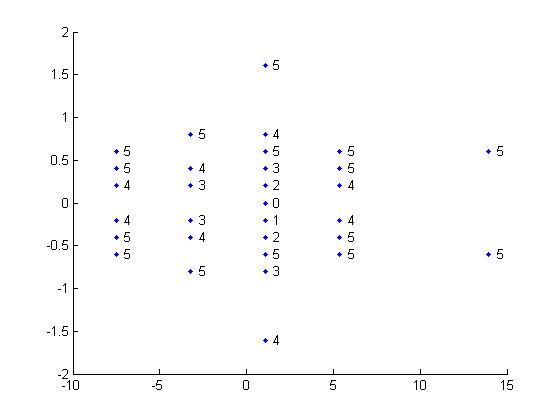} 
   	\caption{First thirty-two  multivalues of multivalued E(k) on complex plane, for b=$\frac{1}{4}$} \label{eagm_e}
    	\vspace{-0.2cm}
	\end{figure}

	The result, for four choices of the sign of the geometric mean (see "fille32" in appendix
	\ref{appendix:program}), is a set of mutivalues whose complex plane map is shown in figure~\ref{eagm_e}.

           We observe that it has the exact same shape as the set of section~\ref{section:K}

	While the shape, in terms of reached integer values of the lattice generating vectors, are the same, the generating vector are of different sizes.

	The horizontal/real factor is obvious: $E(k)/K(k)$, because the real-axis lattice generating vector is $4E(k)$ instead of $4 K(k)$.

\begin{figure}
	\includegraphics[width=0.92\textwidth]{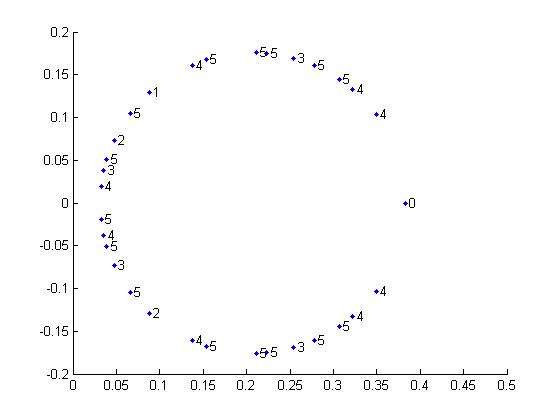} 
   	\caption{First thirty-two multivalues of the multivalued MAGM function $N(b^2)$ on complex plane, for b=$\frac{1}{4}$. The multivalues are located on a circle with midpoint at the x-axis, intersecting the x-axis at $1-N(k^2)$ and $N(b^2)$.}  \label{eagm_n}
    	\vspace{-0.2cm}
	\end{figure}

	The vertical/imaginary factor is $(K(b)-E(b))/K(b)$, because the imaginary-axis lattice generating vector is $4iK(b)-4iE(b)$ instead of $4iK(b)$. The proof of this ratio can be obtained 
	by using the development of the values of the AGM for alternate sign choices, using the same methods used to develop K in the subsections of section~\ref{section:K}, and  
	similar known formulae~\cite{dlmf}.
	The consequence of the unequal scaling is that the ratio between $E(k)$ and $K(k)$, called $N(b^2)=N(1-k^2)$, itself will have multivalues, and those multivalues 
	will be located on a circle on the complex plane. The plot of $N(b^2)$ for 32 multivalues runs in "filln32" program and yields figure~\ref{eagm_n}.
	It is remarkable that $N(b^2)$ forms a dense set even though it is a set built from the ratio of two complex numbers each having 
	a lattice extending to infinity in all directions.

	In fact, it is interesting to note that for any set of points on the complex plane, taking the ratio point-wise between such a scaled version of this set and the original set will yield 
	a circle, if the scaling is different in the real and imaginary directions.

   \vspace{-2mm}
%
%
   \section{MULTIVALUE LATTICES RELATING TO THE INCOMPLETE INTEGRAL OF THE SECOND KIND}
   \label{section:Z}

	The incomplete elliptic integral of the second kind is usually calculated splitting it into a trivial part, proportional to the incomplete elliptic of the first kind, plus a function called Zeta:

	$$
	E(\phi,k) = F(\phi,k)E(k)/K(k) + Z(\phi,k)
	$$

	Just like the computation of $E(k)$, the computation of $Z(\phi,k)$ involves a sum of terms, each one multiplied by a successively larger power of two.

\begin{figure}
	\includegraphics[width=0.7\textwidth]{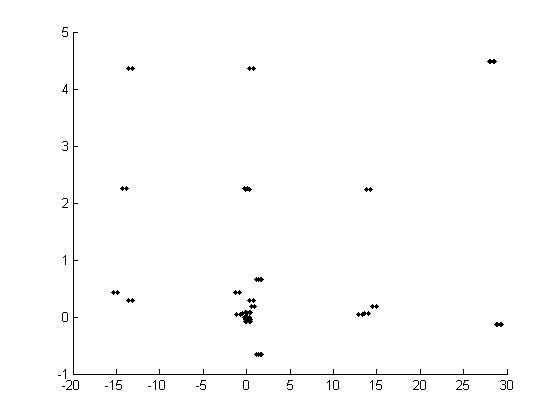} 
   	\caption{First 64 (4 by 4 by 4) multivalues of multivalued $Z(\phi,k)$ on complex plane, for b=$\frac{1}{4}$ and $\sin{\phi}=0.8$} \label{eagm_z_bad}
    	\vspace{-0.2cm}
	\end{figure}

	When we calculate the effect of taking any choice of sign $\sigma_n$, $\delta_n$ and $\gamma_n$, a very chaotic pattern emerges for Z on the complex plane. See figure~\ref{eagm_z_bad}.

	However, if we restrict the choice of sign to $\delta_n$,and make the sign $\gamma$ correlated to $\delta$, with $\gamma_n=\delta_{n-1}$, we find that the pattern of 
	multivalues has regularity - it forms a one-dimensional lattice, see the other figure~\ref{eagm_z_good}.

          The existence of this pattern reflects a property of the Z function, namely that after crossing two branch cuts, it comes back to the same value as before, except for the 
	addition of an imaginary constant. In the following paragraphs, we derive this imaginary constant:

	\begin{figure}
	\includegraphics[width=0.7\textwidth]{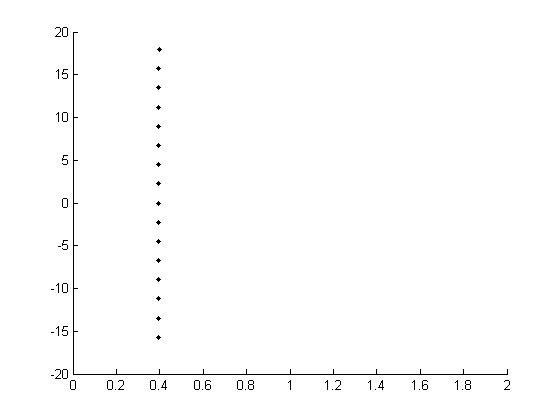} 
   	\caption{First 16 restricted multivalues of multivalued $Z(sin(\phi),k)$ on complex plane,, for b=$\frac{1}{4}$ and $\sin{\phi}=0.8$} \label{eagm_z_good}
    	\vspace{-0.2cm}
	\end{figure}

\subsection {Calculation of the multi-valuedness of Jacobi zeta function.}

	For ease, we will denote by $b$ the function $\sqrt{1-k^2}$ of the elliptic integral parameter $k$, and we may also sometimes use $m$ to denote $k^2$.

	The Jacobi zeta function is multivalued, and the values are separated by a "lattice constant" which can be evaluated. An inverse function of the Jacobi Zeta function would 
	have this "lattice constant" as periodicity. In order to calculate the separation between multiple values of this function, we use its definition as 
	$$ E(\phi,k)-E(k)F(\phi,k)/K(k) $$
	and calculate the change of the function as it crosses a first branch cut, going around the branch point at $\pi/2 + i \sinh^{-1}(b/k)$ or $\pi/2 + i \cosh^{-1}(1/k)$, 
	followed a second branch cut, going around the branch point at the complex conjugate position $\pi/2 - i \sinh^{-1}(b/k)$

	We use the~\cite{FBrCut} reference to measure the change of the function $F$ at the branch cut; this change consists of the addition of this constant $C_F$, along with a switch 
	to the $-F$ function, so as to maintain continuity. When going around the other branch cut, the change consists of the addition of the same constant, along with a switch back 
	to the $F$ function.

	$$  C_F = 4 K(k)-\frac{2}{k}K\left(\frac{1}{k}\right) $$ 

	Similarly, we use the~\cite{EBrCut} reference evaluate the change of the function $E$ after going around the two branch cuts, it is two times this other constant:

	$$ C_E = 4 E(k)-2 k \left(E\left(\frac{1}{k}\right)+\left(\frac{1}{k^2}-1\right)K\left(\frac{1}{k}\right)\right) $$

	By definition of $Z$, the effect on Jacobi $Z$ of going around the two branch points is the addition of a constant $Q_Z = 2C_E - 2 C_F E(k)/K(k)$. One immediately notices 
	that the $4E(k)$ term cancels the scaled $4K(k)$ term. Then one is left with...

	$$Q_Z = 2C_E - 2 C_F \frac{E(k)}{K(k)} = -4k E\left(\frac{1}{k}\right) + 4k K\left(\frac{1}{k}\right) -\frac{4}{k} K\left(\frac{1}{k}\right) + \frac{4}{k} \frac{E(k)K(1/k)}{K(k)}$$

	$$ = -4k E\left(\frac{1}{k}\right) + \left(4k-\frac{4}{k}+\frac{4}{k} \frac{E(k)}{K(k)}\right) K\left(\frac{1}{k}\right) $$

	We further use~\cite{dlmf} to develop the values of the elliptic integral with parameter greater than one into real and imaginary parts, yielding this expression 
	for $Q_Z$:
	
	$$-4 E(k) \mp 4 i E(b) + 4(1-k^2) K(k) \pm 4i k^2 K(b) + \left(4 k^2 -4 +4\frac{E(k)}{K(k)}\right) K(k) \mp i \left(4 k^2 -4 + 4\frac{E(k)}{K(k)}\right) K(b)$$
	
	The real parts of $Q_Z$ cancel out, and it leaves out:
	$$ Q_Z = \mp 4 i E(b) \pm 4i k^2 K(b) \mp  4 i k^2 K(b) \pm 4 i K(b) \mp 4 i E(k)K(b)/K(k) $$ 
	or, removing canceling terms,
	$$ Q_Z = \mp 4i E(b) \pm 4i K(b) \mp 4i E(k)K(b)/K(k) $$
	or, grouping the denominator,
	$$ Q_Z = \mp 4i (E(b)K(k)-K(b)K(k)+E(k)K(b))/K(k)$$
	and so finally, using the Legendre's relation
	$$ Q_Z = \mp 4i \frac{\pi}{2} \Big/ K(k) = \mp \frac{2 \pi i}{K(k)} $$ 

\subsection {Significance of the multivaluedness of Zeta}

	Since the multivaluedness (of the 2D lattice kind) of F is linked to the (double) periodicity of its inverse function, sn u, one can wonder
	 if the multivaluedness of Z in a regular pattern is indicative
	of a periodicity of its inverse function. It is relatively easy to see that it is not the case, at least in the most basic meaning of periodicity. If we take the example 
	$Z(0.5,\sqrt{0.9375}) \approx 0.2920$, meaning the inverse function $\xi$ has a value $\xi(0.2920,\sqrt{0.9375}) \approx 0.5$, it is easy to solve numerically and see that
	$\xi(0.2920+2.2430i,\sqrt{0.9375}) \approx 0.3150+2.7672i \neq 0.5$ because $Z(0.3150+2.7672i,\sqrt{0.9375}) \approx 0.2920+2.2430i$. So $2.2430i$, the unit vector of 
	the 1D lattice of multivalues, is not a period of $\xi$ the inverse function of $Z$.

	However, just as we were able to analytically continue $Z$ around a branch cut, its inverse function $\xi$ could be analytically continued following $Z$ around a branch cut. 
	After following it around two branch cuts (i.e. an even
 	number of times, to restore the original sign), we would arrive to an offset version of $\xi$, which reaches the target value of 0.5 at the precise offset corresponding to the "period". We
	have a kind of periodicity, but it is only when considering $\xi$ and all of its "branches".  


   \vspace{-2mm}
   \section{CONCLUSIONS}
   \label{section:conclu}

	We have shown that using the sign choice ambiguity in AGM calculation yields multivalues that hint at the multivaluedness of the functions we call elliptic integrals.
	This multivaluedness spans a lattice on the complex plane which itself points at the double periodicity of the inverse of these functions, the Jacobian elliptic functions.
	In the case of the incomplete elliptic integral of the second kind, we can only give meaning to some of the multivalues, so many challenges remain.

\newpage
\noindent \textbf{\large Appendices}

\appendix

   \section{A program to calculate multivalues}
   \label{appendix:program}

   \lstset{language = Matlab, 
           basicstyle = \footnotesize}
   \begin{lstlisting}
function [zd,f,kn,ko]=eagm(k,sinphi,signb,sign1s,sign2s,sign3s)
    n0=1;
    g=signb*sqrt(1-k*k);
    a=1;
    u=1/sinphi;
    v=signb*sqrt(1-k*k*sinphi*sinphi)/sinphi;
    zd=0;
    p2=1;
    for j=0:19
        n0=n0-(a*a-g*g)*p2/2;
        if (real(sqrt(z*z-a*a)/z)>=0)
            sig_z=1;
        else
            sig_z=-1;
        end
        if (bitand(2^j,sign3s))
            zd=zd-sig_z*(z-v)/z*sqrt(z*z-a*a)*p2;
        else
            zd=zd+sig_z*(z-v)/z*sqrt(z*z-a*a)*p2;
        end
        t=sqrt(a*g);
        if real(t/(a+g))>=0
            sig_g=1;
        else
            sig_g=-1;
        end
        if bitand(2^j,sign1s) % this integerbit code selects far vs near 
            t=-sig_g*sqrt(a*g); % the far geometric mean
        else
            t=sig_g*sqrt(a*g);  % the near geometric mean
        end
        
        vo=v;
        v=sqrt((u+vo)*(u+vo)-(a-g)*(a-g))/2;
        if real(v/(u+vo))>=0   % sig*v is the root that allows convergence
            sig=1;
        else
            sig=-1;
        end
        if (bitand(2^j,sign2s))
            v=-v*sig;
        else
            v=+v*sig;
        end
        u=u/2+vo/2;
        a=a/2+g/2;
        g=t;
        p2=p2*2;
    end
    f=asin(a/u)/a;
    ko=pi/2/a;
    kn=ko*n0;
end
% fillk2x32
function fillk2x32(k)
if ~exist('k')
    k=sqrt(0.9375);
end
colormap='kbgrmc';
hold on
for jj=31:-1:0
for ll=0:-1:0
for mm=0:-1:0
    if (jj~=16)
        [zd,ff,ee,kc]=eagm(sqrt(0.9375),0.5,1,jj,ll,mm);
        generation=1+floor(log(jj+0.99)/log(2));
        plot(real(kc),imag(kc),strcat(colormap(1+generation),'h'));
        text(real(kc)+0.8,imag(kc),num2str(generation));
    else
        disp(strcat('one point is skipped because it suffers from numeric precision',
                    ' limitations and it would change too much the scale of the plot'));
    end
    [zd,ff,ee,kc]=eagm(sqrt(0.9375),0.5,-1,jj,ll,mm);
    plot(real(kc),imag(kc),'k.');
end
end
end
hold off
end
% fillf8x16
function fillf8x16(sinphi,k)
if ~exist('k')
    k=sqrt(0.9375);
end
if ~exist('sinphi')
    sinphi=0.8;
end
hold on
for jj=7:-1:0
for ll=15:-1:0
for mm=0:-1:0
    [zd,ff,ee,kc]=eagm(k,sinphi,1,jj,ll,mm);
    plot(real(ff),imag(ff),'k.')
end
end
end
hold off
end
% fille32
function fille32(k)
if ~exist('k')
    k=sqrt(0.9375); % example value by defaut
end
hold on
for jj=31:-1:0
for ll=0:-1:0
for mm=0:-1:0
    if (jj~=16)
        [zd,ff,ee,kc]=eagm(k,0.5,1,jj,ll,mm);
        generation=1+floor(log(jj+0.99)/log(2));
        plot(real(ee),imag(ee),'b.');
        text(real(ee)+0.5,imag(ee),num2str(generation));
    else
        disp(strcat('one point is skipped because it suffers from numeric precision',
                    ' limitations and it would change too much the scale of the plot'));
    end
end
end
end
hold off
end
% filln32
function filln32(k)
if ~exist('k')
    k=sqrt(0.9375); % example value by defaut
end
hold on
for jj=31:-1:0
for ll=0:-1:0
for mm=0:-1:0
    if (jj~=16)
        [zd,ff,ee,kc]=eagm(k,0.5,1,jj,ll,mm);
        generation=1+floor(log(jj+0.99)/log(2));
        plot(real(ee/kc),imag(ee/kc),'b.');
        text(real(ee/kc)+0.005,imag(ee/kc),num2str(generation));
    else
        disp(strcat('one point is skipped because it suffers from numeric precision',
                    ' limitations and it would change too much the scale of the plot'));
    end
end
end
end
axis([0 0.5 -0.2 0.2])
hold off
end
% fillz4x4x4
function fillz4x4x4(sinphi,k)
if ~exist('k')
    k=sqrt(0.9375);
end
if ~exist('sinphi')
    sinphi=0.8;
end
hold on
for jj=3:-1:0
for ll=3:-1:0
for mm=3:-1:0
    [zd,ff,ee,kc]=eagm(sqrt(0.9375),0.5,1,jj,ll,mm);
    plot(real(zd),imag(zd),'k.')
end
end
end
hold off
end
%fillz1x16
function fillz1x16(k,sinphi)
if ~exist('k')
    k=sqrt(0.9375);
end
if ~exist('sinphi')
    sinphi=0.8;
end
hold on
for jj=0:-1:0
for ll=15:-1:0
for mm=0:-1:0
    [zd,ff,ee,kc]=eagm2(k,sinphi,1,jj,kk,ll,ll*2);
    plot(real(zd),imag(zd),'k.')
end
end
end
hold off
axis([0 2 -20 20])
end

   \end{lstlisting}
   \vspace{-3mm}
   
   \section{Relationship with the MAGM}
	In this appendix, we will explore another algorithm involving repeatedly taking a square root: the MAGM algorithm discovered by S. Adlaj~\cite{MAGM}.
	First, considering only positive choices of the square root, we prove that that algorithm is equivalent to the algorithm discovered by Gauss (and exposed very clearly
	in many review articles, such as Jameson~\cite{jameson}, and used in our EAGM code). Then we will see this method only works for this choice of square root signs. 

\section{Demonstration of the equivalence of the MAGM with the Gauss method}

	Starting with
	$$a_{n+1}=\frac{a_n+g_n}{2}$$
	and
	$$g_{n+1}=\sqrt{a_n g_n}$$
	we note that
	$$\frac{E(\sqrt{1-b^2})}{K(\sqrt{1-b^2})}=1-S$$
	where
	$$S=\sum_{n=0}^{\infty} 2^{n-1} (a_n^2-g_n^2)$$

	It is easy to see that the first few rows of MAGM calculation correspond exactly to the partial evaluations of S by directly evaluating them explicitly.

\begin{table}[h]
  \footnotesize
  \begin{center}
    \caption{First rows of AGM and MAGM evaluation.}
    \label{tab:table1}
    \begin{tabular}{|l|c|c|c|c|c|c|}
      \hline
      \textbf{row}   & \textbf{arithmetic} & \textbf{geometric} & \textbf{current} & \textbf{MAGM} & \textbf{MAGM} & \textbf{MAGM}\\
      \textbf{index} & \textbf{mean}       & \textbf{mean}      & \textbf{series} &  \textbf{x}    & \textbf{y}    & \textbf{z}   \\
      \hline
      $i$ & $a_i$ & $g_i$ & $1-S$ & $x_i$ & $y_i$ & $z_i$ \\
      \hline
      \hline
      0 & 1               & $b$        & 1                 & 1                 & $b^2$   & 0 \\
      1 & $\frac{1+b}{2}$ & $\sqrt{b}$ & $\frac{1}{2}+\frac{b^2}{2}$ & $\frac{1}{2}+\frac{b^2}{2}$ & $b$ & $-b$  \\
      2 & $\frac{1+2\sqrt{b}+b}{4}$ & $\sqrt{\frac{1+b}{2}\sqrt{b}}$ & $\frac{(1+b)^2}{4}$ & $\frac{(1+b)^2}{4}$ & $\sqrt{b}(1+b)-b$ & $-\sqrt{b}(1+b)-b$ \\
      \ldots & \ldots & \ldots & \ldots & \ldots & \ldots & \ldots \\
     $n+1$ & $\frac{a_n+g_n}{2}$ & $\sqrt{a_n g_n}$ & \ldots & $\frac{x_n+y_n}{2}$ & $z_n+\sqrt{(x_n-z_n)(y_n-z_n)}$ & $2z_n-y_{n+1}$ \\
      \hline
    \end{tabular}
  \end{center}
\end{table}

	In order to establish that this correspondence continues indefinitely, it is possible to make an {\it ad absurdum} demonstration: we start by assuming
	that the correspondence stops being true at some row B, and use the properties of the AGM and MAGM to deduce that the correspondence will stop being true at row B-1.

	First, note that for the correspondence between the series and the MAGM to continue up to row B, all we need is for $x_{n+1}-x_n$ to be equal 
	to $2^{n-1}(a_n^2-g_n^2)$ for $n \leq B$.

	Then note that the scaling invariance of the AGM means we can divide every $(a_n , g_n )$ pair by $\frac{1+b}{2}$ 
	and still have an AGM sequence of pairs of numbers.

	The translational invarience of the MAGM means we can add b to every triplet, the resulting sequence of triplets still being MAGM.

	And note that the scaling invariance of the MAGM means we can multiply every $(x_n, y_n, z_n)$ triplet by $\frac{(1+b)^2}{2}$ and still have a MAGM 
	sequence of triplets.

	If we change the origin of indexing of the AGM by one, then the consecutive terms of $S$ each change by a factor of two, because of the exponent of 2 
	in each term $2^{n-1} (a_n^2-g_n^2)$ changes by one.

	Combining the scalings, translation, and index origin change by one, we note that the second row effectively becomes a first row, only with the 
	symbol $b$ replaced by a symbol $v$, where $v=\frac{2 \sqrt{b}}{1+b}$.

	We get the impossible conclusion that the correspondence stops being true at $B-1$ as far as $v$ is considered, whereas it stopped being true 
	at $B$ as far as $b$ is considered. Therefore, this correspondence must continue infinitely.

\section{What about $g_{n+1}=- \sqrt{a_n g_n}$ at specific steps in the calculation}

	We ran the MAGM, letting the sign of the square root of $(x_n-z_n)(y_n-z_n)$ be negative for various n (combined with many choice combinations of $\sigma_n$), and we did not obtain any good multivalues of E(k).

	Thinking about the above demonstration, where we transform a second row into a first row, such a transformation is only possible if all the sign choices for the square root are the same, positive or negative.

	All the signs being negative never yields convergence, so it is no wonder that the MAGM algorithm only delivers E(k) for the "all positive sign choices" case.

\bibliographystyle{amsplain}

\end{document}